\documentclass[12pt]{amsart}
\usepackage{amsthm}
\swapnumbers

\makeatletter
\def\swappedhead#1#2#3{%
  \thmnumber{\@upn{\the\thm@headfont#2\@ifnotempty{#1}{.~}}}%
  \thmname{#1}%
  \thmnote{ {\the\thm@notefont(#3)}}}
\makeatother

\newtheoremstyle{dotless-thm}
  {3pt}
  {3pt}
  {\itshape}
  {}
  {\bfseries}
  {}
  {.5em}
  {}
\theoremstyle{dotless-thm}

\usepackage[margin=1in]{geometry}          
\usepackage{amsmath, amscd, amssymb}     
\usepackage{amsfonts, lipsum}    
\usepackage{amsthm}
\usepackage[all]{xy}               

\numberwithin{equation}{section}

\def\m{\medskip}

\newtheorem{thm}{Theorem}[section]

\newtheorem{lem}[thm]{Lemma}
\newtheorem{prop}[thm]{Proposition}
\newtheorem{cor}[thm]{Corollary}

\newtheorem{remark}[thm]{Remark}
\newtheorem{remarks}[thm]{Remarks}

\newtheorem{quest}[thm]{Question}

\newtheorem{df}[thm]{Definition}

\newcommand\theoref{Theorem~\ref}
\newcommand\lemref{Lemma~\ref}
\newcommand\propref{Proposition~\ref}
\newcommand\corref{Corollary~\ref}

\def\SB{\operatorname{SB}}
\def\cat{\operatorname{cat}}
\def\ballcat{\operatorname{ballcat}}

\def\crit{\operatorname{Crit}}
\def\pt{\operatorname{pt}}
\def\cl{\operatorname{cl}}

\def\gl{\lambda}
\def\A{\mathcal A}
\def\N{\mathbb N}

\def\R{\mathbb R}
\def\Z{\mathbb Z}

\def\ts{\times}
\def\wt{\widetilde}
\def\p{{\bf Proof. }}

\long\def\forget#1\forgotten{} %

\begin{document}

\title[Maps of Degree 1 and LS Category]{Maps of Degree 1 and Lusternik--Schnirelmann Category }

\author{Yuli B. Rudyak} 
\address{Department of Mathematics\newline
University of Florida\newline
P.O. Box 118105\newline
Gainesville, FL 32611-8105\newline
 USA}
 
 \email{rudyak@ufl.edu}

\begin{abstract}
Given a map $f: M \to N$ of degree 1 of closed manifolds, is it true that the  Lusternik--Schnirelmann category of the range of the map is not more that the category of the domain? We discuss this and some related questions.
\end{abstract}

\maketitle

\section{Introduction}

 Below $\cat$ denotes the (normalized) Lusternik--Schnirelmann category,~\cite{CLOT}.
 
 \m
\begin{quest}\label{question}
Let $M, N$, $\dim M=\dim N=n$, be two closed connected orientable manifolds, and let $f: M \to N$ be a map of degree $\pm 1$. Is it true that $\cat M \geq \cat N$? 
\end{quest}

There are several reasons to conjecture that the above-mentioned question has the affirmative answer for all maps $f$ of degree 1. Indeed, informally, the domain of $f$ is more ``massive'' than the range of $f$. For example, it is well known that the induced maps $f_*: H_*(M) \to H_*(N)$ and  $f_*: \pi_1(M) \to \pi_1(N)$ are surjective. 

\m In~\cite{R2} I proved some results that confirm the conjecture under some suitable hypotheses, and now people speak about the Rudyak conjecture, cf. \cite[Open Problem 2.48]{CLOT},~ \cite{D}, (although I prefer to state questions rather than conjectures). To date we do not know any counterexample.

\m Here we demonstrate more situations when the conjecture holds, 
and discuss some variations of the conjecture.

\section{Preliminaries}

\begin{df}\rm Let $M,N, \dim M=\dim N=n$, be two closed connected oriented smooth manifolds, and let $[M]\in H_n(M)=\Z$ and $[N]\in H_n(N)=\Z$ be the corresponding fundamental classes. Given a map $f: M\to N$, we define the {\em degree} of a map $f:M \to N$ as the number $\deg f\in \Z$ such that $f_*[M]=(\deg f) [N]$. 
\end{df}

Let $E$ be a commutative ring spectrum (cohomology theory).

\begin{lem}\label{mono} If $M$ is $E$-orientable and $f: M \to N$ is a map of degree $1$ then $f_*\colon E_*(M)\to E_*(N)$ is an epimorphism and $f^*\colon E^*(N)\to E^*(M)$ is a monomorphism.  
\end{lem}

\p See~\cite[Theorem V.2.13]{R1}, cf. also Dyer~\cite{Dyer}.
\qed

\begin{remark}\label{r:local}\rm An analog of \lemref{mono} holds for cohomology with local coefficients.  Let $\A$ be a local coefficient system of abelian groups on $N$, and let $f^*(\A)$ be the induced coefficient system. Then $f_*:H_*(M;f^*\A)\to H_*(N;\A)$ is a split epimorphism and $f^*\colon H^*(N;\A)\to H^*(M;f^*\A)$ is a split monomorphism. The proofs are based on Poincar\'e duality with local coefficients and the equality $f_*(f^*x\frown y)=x\frown f_* y$ for  $x\in H^*(N;\A), y\in H_*(M;f^*\A)$. See e.g.~\cite{B}.
\end{remark}

\begin{df}\rm  Given a CW space $X$, the {\em Lusternik--Schnirelmann category} $\cat X$ of $X$ is the least integer $m$ such that there exists an open covering $\{A_0, A_1,\ldots, A_m \}$ with each $A_i$ contractible in $X$ (not necessary in itself). If no such covering exist, we put  $\cat X=\infty$. 
\end{df}

 Note the inequality $\cat X \leq \dim X$ for $X$ connected.

\m In future we abbreviate Lusternik--Schnirelmann to LS. A good source for LS theory is~\cite{CLOT}.

\m Given a closed smooth manifold $M$ and a smooth function $f: M \to \R$, the number of critical points of $f$ can't be less than $\cat M$,~\cite{LS1, LS2}. This result turned out to be the starting point of LS theory. Currently, the LS theory is a broad area of intensive topological research.

\m Let $X$ be a path connected space. Take a point $x_0\in X$, set 
$$
PX=P(X,x_0)=\{\omega\in X^{[0,1]} \bigm| \omega(0)=x_0\}
$$
and consider the fibration $p: PX \to X, p(\omega)=\omega(1)$. Let $p_k=p_{k}^X: P_k(X)\to X$ be the $k$-fold fiberwise join $PX*_X\cdots *_X* PX \to X$. According to the Ganea--\v Svarc theorem, \cite{CLOT, S},  the inequality $\cat X<k$ holds if and only if the fibration $p_k: P_k(X)\to X$ has a section. In other words, the number $\cat X$ is the least $k$ such that the fibration $p_{k+1}: P_{k+1}(X)\to X$ has a section.

\section{Approximations}

Recall (see the Introduction) that we discuss whether the existence of a map $f: M \to N$ of degree 1 implies the inequality $\cat M \geq \cat N$. Here we present two  results appearing when we approximate the LS category by the cup-length and Toomer invariant. 

\m Recall the following cup-length estimate of LS category. Let $E$ be a commutative ring spectrum (cohomology theory). The cup-length of $X$ with respect to $E$ is the number 
\[
\cl_E(X):=\sup \{m\bigm| u_1\smallsmile\cdots \smallsmile u_m \neq 0 \text{ where } u_i\in\wt E^*(X)\}.
\]
The well-known cup-length theorem~\cite{CLOT} asserts that $\cl_E(X)\leq \cat X$.

\m We give another estimate of  LS category. 

\begin{df}\label{d:toomer} \rm Define the (cohomological) Toomer invariant 
\[
e_E^*(X)=\sup\{k\bigm | \ker \{p_k^*: E^*(X))\to E^*(P_k(X))\}\neq 0\}.
\] 
\end{df}

Note the decreasing sequence $\cdots \supset \ker (p_k^*)\supset \ker (p_{k+1}^*)
\supset \cdots$. Moreover,
\[
e_E^*(X)=\inf\{k\bigm | \ker \{p_k^*: E^*(X)\to E^*(P_k(X))\}= 0\}-1.
\] 

\m In the definition of the Toomer invariant, $E$ does not need to be a ring spectrum, it can be an arbitrary spectrum.
 
\begin{prop}\label{p:doubleineq}
 We have $\cl_E(X)\leq e_E^*(X)\leq \cat X$.
\end{prop}

\p First, $p_m=p_{m}^X$ has a section for $m> \cat X$, and so $\ker p_{m}^*=0$ for all $m>\cat X$. Hence, $e_E^*(X)\leq \cat X$. Now, we put $\cl_E(X)=l$, $e_E^*(X)=k$ and prove that $l\leq k$. Take $u_1, \ldots, u_l\in \wt E^*(X)$ such that $u_1\smallsmile \cdots \smallsmile u_l\neq 0$. Then, since $p_{k+1}^*$ is monic, we conclude that $p_{k+1}^*(u_1\smallsmile \cdots \smallsmile u_l)\neq 0$. In other words, 
\[(p_{k+1}^*u_1)\smallsmile \cdots \smallsmile  (p_{k+1}^* u_l)\neq 0
\]
 in $P_{k+1}X$. Hence, $\cat P_{k+1}(X)\geq l$ because of the cup-length theorem. It remains to note that $\cat P_{k+1}(X)\leq k$ for all $k$, see~\cite[Prop. 1.5(ii)]{R2} or~\cite{CLOT}.
 \qed

\begin{prop}\label{p:ineq} Let $M$ be a closed connected $E$-orientable manifold, and let $f: M \to N$ be a map of degree $\pm 1$. Then $\cl_E(M)\geq \cl_E(N)$ and $e^*_E(M)\geq e^*_E(N)$.
\end{prop}

\p Put $\cl_E(N)=l$ and take $u_1, \ldots u_l\in \wt E^*(N)$ such that $u_1\smallsmile \cdots \smallsmile u_l\neq 0$. Then $f^*(u_1\smallsmile \cdots \smallsmile u_l)\neq 0$ by \ref{mono}. So, $f^*(u_1)\smallsmile \cdots \smallsmile f^*(u_l)\neq 0$ in $M$, and thus $\cl_E(M)\geq l$.

\m Now, let $e^*_E(M)=k$. Then $(p_i^M)^*$ is a monomorphism for $i>k$. Consider the commutative  diagram
\[
\CD
P_iM @>P_if>> P_i N\\
@Vp_i^M VV @Vp_i^N VV\\
M@>f>> N.
\endCD
\]
and note that, by \ref{mono}, $f^*((p_i^M)^*)$ is monic for $i>k$. Hence, because of commutativity of the diagram, $(p_i^N)^*$ is monic for $i>k$. Thus, $e^*_E(N)\leq k$. 
\qed

\begin{cor}\label{c:cat=cl}
Let $M$ be a closed connected $E$-orientable manifold, and let $f: M \to N$ be a map of degree $\pm 1$.

{\rm(i)} Suppose that $e^*_E(N)=\cat N$. Then $\cat M\geq \cat N$.

 {\rm(ii)} Suppose that $\cl_E(N)=\cat N$. Then $\cat M\geq \cat N$. 
\end{cor}

\p (i) We have $\cat M \geq e^*_E(M)\geq e^*_E(N)=\cat N$, the second inequality following from \propref{p:ineq}. 

(ii) Because of (i), it suffices to prove that  $e^*_E(N)=\cat N$. But this holds since $\cl_E(N)\leq e^*_E(N)\leq \cat N$.
\qed

\section{Low-Dimensional Manifolds}\label{s:low} 

We prove that for $n\leq 4, \cat M^n \geq \cat N^n$ provided that there exists a map $f: M \to N$ of degree 1. The inequality holds trivially for $n=1$. 

\m The case $n=2$ is also simple. Denote by $g(X)$ the genus of a closed connected orientable surface $X$. Then $g(M)\geq g(N)$ because of surjectivity of $f_*: H_2(M) \to  H_2(N)$, see~\ref{mono}.  Furthermore, $\cat X=1$ if $g(X)=0$ ($ X=S^2)$ and $\cat X=2$ for $g>1$. 

\m The case $n=3$ is considered in~\cite[Corollary 1.3]{OR}.

\m The case $n=4$.
First, if $\cat N=4$ then $\cat N=\cat M$ by~\cite[Corollary 3.6(ii)]{R2}. 

\m
Next we consider $\cat N=3$ and prove that $\cat M \geq 3$. By way of contradiction, assume that $\cat M\leq 2$, and hence the group $\pi_1(M)$ is free by~\cite{DKR}. 
Now, $\pi_1(N)$ is free since $\deg f=1$, see~\cite{DR}. But then $N$ has a CW decomposition whose 3-skeleton is a wedge of spheres,~\cite{Hil}, and hence $\cat N\leq 2$, a contradiction. Finally, the case $\cat M=1<2=\cat N$ is impossible for trivial reasons ($M$ should be a homotopy sphere, and therefore $N$ should be a homotopy sphere).

\section{Some exemplifications}

The following result is a weak version of \cite[Theorem 3.6(i)]{R2}.

\begin{thm} \label{t:main}
Let $M, N$ be two smooth closed connected stably parallelizable manifolds, and assume that there exists a map $f: M \to N$ of degree $\pm 1$. If $N$ is $(q-1)$-connected and $\dim N\leq 2q\cat N-4$ then $\cat M \geq \cat N$.
\qed
\end{thm}

\begin{remark}\rm
 In~\cite[Theorem 3.6(i)]{R2} we use~\cite[Corollary 3.3(i)]{R2} where, in turn, we require $\dim N \geq 4$. However, the case $\dim N \leq 4$ is covered by Section~\ref{s:low}.
 \end{remark}

\begin{thm} \label{t:torus}
Let $T^k$ denote the $k$-dimensional torus. Let $M, N$ be two smooth closed connected stably parallelizable  manifolds, and assume that there exists a map $f: M \to N$ of degree $\pm 1$. Then there exists $k$ such that $\cat (M \ts T^k)\geq \cat (N\ts T^k)$. 
\end{thm}

\p Put $\dim M=n$ and note that $\cat (T^k\ts M)\geq \cat T^k=k$. Now, if $k\geq n+4$ then
\[
2\cat M-4\geq 2k-4 \geq k+n,
\]
and we are done by \theoref{t:main}.
\qed

\m Another example. Consider the exceptional Lie group $G_2$. Recall that $\dim G_2=14$. Note that  $G_2$ is parallelizable being a Lie group.

\begin{prop} Let $M$ be a stably parallelizable 14-dimensional closed manifold that admits a map $f: M \to G_2$ of degree $\pm 1$. Then $\cat M \geq \cat G_2$.
\end{prop}

\p The group $G_2$ is 2-connected and $\cat G_2=4$,~\cite{IM}. Now the result follows from \theoref{t:main} with $N=G_2$ and $q=3$ because $14=\dim G_2\leq 2q\cat G_2-4=20$.
\qed

\m Let $SO_n$ denote the special orthogonal group, i.e., the group of the orthogonal $n\times n$-matrices of determinant 1. Recall that  $\dim SO_n=n(n-1)/2$.

\begin{thm} \label{t:son}
Let $M$ be a closed connected smooth  manifold, and let $f: M \to SO_n$ be a map of degree $\pm 1$. Then $\cat M \geq \cat(SO_n)$ for $n\leq 9$.
\end{thm}

\p We apply \corref{c:cat=cl} for the case $E=H\Z/2$, i.e., to arbitrary closed connected manifolds.  Below $\cl$ denotes $\cl_{\Z/2}$. Because of \corref{c:cat=cl} it suffices to prove that $\cat SO_n=\cl(SO_n)$ for $n\leq 9$. Recall that $H^*(SO_n;\Z/2)$ is the polynomial algebra on generators $b_i$ of odd degree $i<n$, truncated by the relations ${b_i^{p_i}}=0$ where $p_i$ is the smallest power of 2 such that ${b_i^{p_i}}$ has degree $\geq n$,~\cite{H}. In other words,
\[
H^*(SO_n;\Z/2)=\Z/2[b_1, \ldots, b_k, \ldots]/(b_1^{p_1},  \ldots, b_k^{p_k}, \ldots).
\]
Note  that $p_k=1$ for $2k-1>n$, and so $H^*(SO_n)$ is really a truncated polynomial ring (not a formal power series ring). Hence, 
 \[
 \cl(SO_n)=(p_1-1)+\ldots +(p_k-1)+\ldots
 \]
 and the sum on the right is finite because $p_k=1$ for all but finitely many $k$'s.
 
\m For sake of simplicity,  we use the notation $S_n$ for  $H^*(SO_n;\Z/2)$. We have:

$\dim SO_3=3$, $\cl(SO_3)=3 \text{ because } S_3=\Z/2[b_1]/(b_1^4)$.\\
$\dim SO_4=6$, $\cl(SO_4)=4 \text{ because } S_4=\Z/2[b_1,b_3]/(b_1^4,  b_3^2)$.\\
$\dim SO_5=10$, $\cl(SO_5)=8 \text{ because } S_5=\Z/2[b_1,b_3]/(b_1^8,  b_3^2)$.\\
$\dim SO_6=15$, $\cl(SO_6)=9 \text{ because } S_6=\Z/2[b_1,b_3, b_5]/(b_1^8,  b_3^2, b_5^2)$.\\
$\dim SO_7=21$, $\cl(SO_7)=11\text{ because } S_7=\Z/2[b_1,b_3, b_5]/(b_1^8,  b_3^4, b_5^2)$.\\
$\dim SO_8=28$, $\cl(SO_8)=12 \text{ because } S_8=\Z/2[b_1,b_3, b_5, b_7]/(b_1^8,  b_3^4, b_5^2, b_7^2)$.\\
$\dim SO_9=36$, $\cl(SO_9)=20 \text{ because } S_9=\Z/2[b_1,b_3, b_5, b_7]/(b_1^{16},  b_3^4, b_5^2, b_7^2)$.\\

The values of $\cat SO_n, n\leq 9$ are calculated in~\cite{IMN}, see also~\cite{I2}. Compare these values with the above-noted values of $\cl(SO_n)$ and see that $\cat SO_n=\cl(SO_n)$ for $n\leq 9$. 
\qed

\begin{remarks}\rm  1. The anonymous referee noticed that, probably, the method of \theoref{t:son} can also be applied to other Lie groups ($U_n, SU_n$, etc.). This is indeed true, but we do not develop these things here. 

 \m 2. In above-mentioned \theoref{t:main} we can relax the assumption on $M$ by requiring that the normal bundle of $M$ be stably fiber homotopy trivial, i.e., that $M$ is $S$-orientable where $S$ is the sphere spectrum. However, we can't provide the same weakening for $N$ because our proof in~\cite{R2} uses surgeries on $N$. 
\end{remarks}

\section{Theme and Variations: Other Numerical Invariants Similar to LS Category}

Let $\crit X$ denote the minimum number of critical points of a smooth function $f:X \to \R$ on the closed smooth manifold $X$, and let $\ballcat (X)$ be the minimum $m\in \N$ such that there is a covering of $X$ by $m+1$ smooth open balls. It is known that $\cat X \leq \ballcat(X)\leq \crit M-1$, see~\cite{CLOT}. Note that there are examples with $\cat M<\ballcat(M)$. Indeed, there are examples of manifolds $M$ such that $\cat M=\cat (M\setminus \pt)$,~\cite{I1}, while $\cat (M\setminus \pt)=\cat M-1$ whenever $\cat M=\ballcat(M)$. On the other hand, there are no known examples with $\ballcat M +1<\crit M$.

 \m Now,  we can pose an open question whether $\ballcat M\geq \ballcat N$ and $\crit M \geq \crit N$ provided there exists a map $f: M \to N$ of degree 1.
 
\m We can also consider the number $\crit^*(X)$, that is, the minimum number of {\em nongenerate} critical points of a smooth function $f:X \to \R$ on the closed smooth manifold $X$. There is a big difference between $\crit X$ and $\crit^*X$. For example, if $S_g$ is a surface of genus $g\geq 1$ then $\crit S_g=3$ while $\crit^*S_g=2g$. So, we can ask if $\crit^* M \geq \crit^* N$ provided there exists a map $f: M \to N$ of degree~1. One of the lower bounds of $\crit^*(X)$ is the sum of Betti numbers $\SB(X)$, the inequality $\SB(X)\leq\crit^*(X)$ being a direct corollary of Morse theory, \cite{M}, and we can regard $\SB(X)$ as an approximation of $\crit^*X$. Now, if $f: M \to N$ is a map of degree~1 then the inequality $\SB(M)\geq \SB(N)$ follows from \ref{mono}.

\m If $M,N$ are closed simply-connected manifolds of dimension $\geq 6$ and there exists a map $M\to N$ of degree~1 then $\crit^*(M)\geq \crit^*(N)$. Indeed, for every Morse function $h: X\to  \R$ on a closed connected smooth manifold $X$ we have the Morse inequalities 
\[
m_{\gl}\geq r_{\gl}+t_{\gl}+t_{\gl-1}
\]
 where $m_{\gl}$ is the number of critical points of index $\gl$ of $h$ and $r_{\gl}, t_\lambda$ are the rank and the torsion rank of $H_{\gl}(X)$, respectively.  Now, if $X$ is a closed simply connected manifold with $\dim M \geq 6 $ then $X$ possesses a Morse function for which the above-mentioned Morse inequalities turn out to be equalities. This is a well-known Smale Theorem~\cite{Sm}. Now, the inequalities $r_{\gl}(M)\geq r_{\gl}(N)$ and $t_{\gl}(M) \geq t_{\gl}(N)$ follow from \ref{r:local}, and we are done.
\qed

\m {\bf Acknowledgments:} The work was partially supported by a grant from the Simons Foundation (\#209424 to Yuli Rudyak). I am very grateful to the anonymous referee who assisted a lot in my work with \corref{c:cat=cl} and \theoref{t:son},  and made other valuable suggestions that improved the paper.

\end{document}